\documentclass[11pt,a4paper]{amsart}
\usepackage{amsfonts,amsmath,amssymb,amsthm}
\usepackage{times}   

\numberwithin{equation}{section}

\DeclareSymbolFont{SY}{U}{psy}{m}{n}
\DeclareMathSymbol{\emptyset}{\mathord}{SY}{'306}

\DeclareMathOperator{\Ran}{Ran}
\DeclareMathOperator{\Dom}{Dom}
\DeclareMathOperator{\spec}{spec}

\newcommand{\norm}[1]{\|#1\|}
\newcommand{\Norm}[1]{\left\|#1\right\|}
\newcommand{\abs}[1]{|#1|}

\newcommand{\R}{\mathbb{R}}
\newcommand{\EE}{\mathsf{E}}
\newcommand{\ii}{\mathrm i}

\newcommand{\cC}{{\mathcal C}}
\newcommand{\cD}{{\mathcal D}}
\newcommand{\cG}{{\mathcal G}}
\newcommand{\cH}{{\mathcal H}}
\newcommand{\cK}{{\mathcal K}}
\newcommand{\cL}{{\mathcal L}}
\newcommand{\cN}{{\mathcal N}}

\newcommand{\cS}{{\mathcal S}}
\newcommand{\cT}{{\mathcal T}}

\newtheorem{theorem}{Theorem}[section]{\bf}{\it}
\newtheorem{corollary}[theorem]{Corollary}{\bf}{\it}
\newtheorem{example}[theorem]{Example}{\it}{\rm}
\newtheorem{lemma}[theorem]{Lemma}{\bf}{\it}
\newtheorem{remark}[theorem]{Remark}{\it}{\rm}
\newtheorem{hypothesis}[theorem]{Hypothesis}{\bf}{\it}

\title[Reducing graph subspaces and strong solutions to Riccati equations]
{Reducing graph subspaces and strong solutions to operator Riccati equations}
\date{\today}

\author[K.\ A.\ Makarov]{Konstantin A.\ Makarov}
\address{K.~A.~Makarov, Department of Mathematics, University of Missouri, Columbia, MO 65211, USA}
\email{makarovk@missouri.edu}

\author[S.\ Schmitz]{Stephan Schmitz$^*$}
\address{S.~Schmitz,
FB 08 - Institut f\"{u}r Mathematik, Johannes Gutenberg-Universit\"{a}t Mainz,
Stau\-dinger Weg 9, D-55099 Mainz, Germany}
\email{schmist@uni-mainz.de}

\author[A.\ Seelmann]{Albrecht Seelmann$^\dagger$}
\address{A.~Seelmann,
FB 08 - Institut f\"{u}r Mathematik, Johannes Gutenberg-Universit\"{a}t Mainz,
Staudinger Weg 9, D-55099 Mainz, Germany}
\email{seelmann@mathematik.uni-mainz.de}

\thanks{$^*,^\dagger$ Individual parts of the  considerations presented in this work will be contained in the author's
Ph.D.\ thesis.}

\begin{document}

\begin{abstract}
The problem of block diagonalization for diagonally dominant symmetric block operator matrices with self-adjoint diagonal entries
is considered. We show that a reasonable block diagonalization with respect to a reducing graph subspace requires a related
skew-symmetric operator to be a strong solution to the associated Riccati equation. Under mild additional regularity conditions, we
also establish that this skew-symmetric operator is a strong solution to the Riccati equation if and only if the graph subspace is
reducing for the given operator matrix. These regularity conditions are shown to be automatically fulfilled whenever the
corresponding relative bound of the off-diagonal part is sufficiently small. This extends the results by Albeverio, Makarov, and
Motovilov in [Canad.\ J.\ Math.\ Vol.\ \textbf{55}, 2003, 449--503], where the off-diagonal part is required to be bounded.
\end{abstract}

\maketitle

\section{Introduction}
In the present work we address the problem of block diagonalization  for possibly unbounded operator matrices in a Hilbert space
$\cH$ of the form
\begin{equation}\label{eq:blockOp}
 B = \begin{pmatrix} A_0 & W\\ W^* & A_1 \end{pmatrix} = \begin{pmatrix} A_0 & 0\\0 & A_1 \end{pmatrix} +
     \begin{pmatrix}0 & W\\ W^* & 0\end{pmatrix}=:A+V
\end{equation}
with respect to a given orthogonal decomposition $\cH=\cH_0\oplus\cH_1$. Here, the diagonal part $A$ of $B$ is a self-adjoint
operator on
\[
 \Dom (A)= \Dom(A_0)\oplus \Dom(A_1)\,,
\]
the off-diagonal part $V$ is a densely defined symmetric operator on a possibly larger domain 
\[
 \Dom(V)\supset \Dom(A)\,,
\]
and $B $ is understood as the diagonally dominant operator sum
\[
 B = A+V\quad \text{on} \quad \Dom(B) = \Dom(A)\,.
\]
For the concept of diagonally dominant operator matrices we refer to \cite[Section 2.2]{Tre08}.

The problem of block diagonalization for operator matrices of the form \eqref{eq:blockOp} is closely related to the existence of
reducing graph subspaces for $B$.

Recall that a closed subspace $\cG\subset\cH$ is said to be \emph{invariant} for a linear operator $B$ if
\[
 \Ran\bigl(B|_{\Dom(B)\cap\cG}\bigr) \subset \cG\,.
\]
If a closed subspace $\cG$ and its orthogonal complement $\cG^\perp$ are invariant for $B$, it is natural to consider the
restrictions $B|_{\Dom(B)\cap\cG}$ and $B|_{\Dom(B)\cap\cG^\perp}$. However, the original operator $B$ coincides with the direct
sum $B|_{\Dom(B)\cap\cG}\oplus  B|_{\Dom(B)\cap\cG^\perp}$ if and only if $\Dom(B)$ splits as
\begin{equation}\label{eq:domSplitting}
 \Dom(B) = \bigl(\Dom(B)\cap\cG\bigr) + \bigl(\Dom(B)\cap\cG^\perp\bigr)\,.
\end{equation}

In this work, following the standard terminology, we say that a closed subspace $\cG\subset\cH$ is \emph{reducing} for $B$ if both
$\cG$ and its orthogonal complement $\cG^\perp$ are invariant for $B$ and the splitting property \eqref{eq:domSplitting} holds.
Note that this splitting  property of the domain is not self-evident if $B$ is unbounded, even if it is self-adjoint
(see \cite[Example 1.8]{Schm12} for a counterexample).

Throughout this work, we are only interested in \emph{graph subspaces} $\cG\subset\cH$ that are associated with bounded operators
from $\cH_0$ to $\cH_1$, that is,
\[
 \cG=\{f\oplus Xf \mid f\in\cH_0\}=:\cG(\cH_0,X)
\]
for some bounded linear operator $X\colon\cH_0\to\cH_1$. For a discussion of a more general concept of graph subspaces where the
``angular operator" $X$ is allowed to be unbounded or even non-closable, we refer to \cite{KM05}.

We establish that the property of a graph subspace $\cG(\cH_0,X)$ to be reducing for the operator $B=A+V$ can be characterized by
the following chain of operator extension relations:
\begin{equation}\label{eq:sandwichIntro}
 (I_\cH-Y)^{-1} (A-YV)(I_\cH- Y)\subset   (A+V) \subset (I_\cH+Y)(A+VY)(I_\cH+Y)^{-1}\,,
\end{equation}
where $I_{\cH}$ denotes the identity operator on $\cH$ and $Y$ is the bounded skew-symmetric operator on $\cH=\cH_0\oplus\cH_1$
given by the $2\times 2$ block operator matrix
\[
 Y = \begin{pmatrix} 0 & -X^*\\ X & 0 \end{pmatrix}\,.
\]
In spite of the fact that the operators $A-YV$ and $A+VY$ are block diagonal with respect to the orthogonal decomposition
$\cH=\cH_0\oplus\cH_1$ (see eqs. \eqref{eq:blockDiag1stIntro} and \eqref{eq:blockDiag2ndIntro} below), the extension relations in
\eqref{eq:sandwichIntro} do not provide satisfying information towards a block diagonalization for the operator $A+V$ on its
natural domain. In fact, in order to have a reasonable block diagonalization for the operator $A+V$, it is natural to require that
at least one of the two relations in \eqref{eq:sandwichIntro} is an operator equality. Moreover, it turns out that if the left-hand
relation in \eqref{eq:sandwichIntro} is an equality, then so is the right-hand one. Thus not all logically possible outcomes in
\eqref{eq:sandwichIntro} may occur (see Theorem \ref{thm:blockDiag} and the discussion in Subsection \ref{subsec:diss}).

One of our principal results states that the chain of operator equalities
\begin{equation*}
 (I_\cH-Y)^{-1} (A-YV)(I_\cH-Y)=  (A+V) = (I_\cH+Y)(A+VY)(I_\cH+Y)^{-1}
\end{equation*}
holds if and only if
\begin{enumerate}
 \renewcommand{\theenumi}{\roman{enumi}}
 \item the graph subspace $\cG(\cH_0,X)$ is reducing for $A+V$,
 
 and

 \item $Y$ is a \emph{strong solution to the operator Riccati equation} $AY-YA-YVY+V=0$, that is,
       $\Ran(Y|_{\Dom(A)})\subset\Dom(A)$ and
       \[
        AYx - YAx - YVYx + Vx = 0 \quad\text{ for }\quad x\in\Dom(A)\,.
       \]
\end{enumerate}
In this case, the operator $A+V$ admits the two block diagonalizations
\begin{equation}\label{eq:blockDiag1stIntro}
 (I_{\cH}+Y)^{-1} (A+V) (I_{\cH}+Y) = A+VY = \begin{pmatrix} A_0+WX & 0\\ 0 & A_1-W^*X^* \end{pmatrix}
\end{equation}
and 
\begin{equation}\label{eq:blockDiag2ndIntro}
 (I_{\cH}-Y)(A+V) (I_{\cH}-Y)^{-1}  = A-YV = \begin{pmatrix} A_0 + X^*W^* & 0\\ 0 & A_1-XW\end{pmatrix}\,.
\end{equation}
In particular, the two block diagonalizations lead to one single unitary block diagonalization with respect to the decomposition
$\cH=\cH_0\oplus\cH_1$,
\begin{equation}\label{eq:uniblockdiag}
 U^* (A+V) U=\begin{pmatrix} B_0 & 0\\ 0 & B_1 \end{pmatrix}\,,
\end{equation}
where $U$ is the unitary transformation from the polar decomposition $I_{\cH}+Y=U\abs{I_{\cH}+Y}$ (hence
$I_{\cH}-Y=\abs{I_{\cH}+Y} U^*$), and the diagonal entries $B_0$ and $B_1$ are similar to the corresponding diagonal entries of the
right-hand sides of \eqref{eq:blockDiag1stIntro} and \eqref{eq:blockDiag2ndIntro} (see Remark \ref{rem:unitary}).

Although we are not able to prove that (i) and (ii) are equivalent in general, we do not have an appropriate counterexample for
which only one of the two requirements holds. Hence, in the framework of our approach, it remains an open problem whether the
requirements (i) and (ii) are (logically) independent or not.

However, under natural additional regularity conditions that automatically hold if, for instance, the off-diagonal perturbation $V$
is bounded, we show that the properties (i) and (ii) are indeed equivalent. In this case, the operator $A+V$ can be diagonalized in
either of the two forms \eqref{eq:blockDiag1stIntro} and \eqref{eq:blockDiag2ndIntro} by making use of the similarity
transformations $(I_\cH+Y)$ and $(I_\cH-Y)^{-1}$, respectively. This extends the results by Albeverio, Makarov, and Motovilov from
\cite[Section 5]{AMM03}, where the operator $V$ is assumed to be bounded and only the block diagonalization
\eqref{eq:blockDiag1stIntro}, but not \eqref{eq:blockDiag2ndIntro}, is mentioned. At the same time, our proof fills in a gap in
reasoning in the proof of \cite[Lemma 5.3]{AMM03} (see the proof of Theorem \ref{thm:AMM} and Remark \ref{rem:AMM}).

The paper is organized as follows.

In Section 2, we show that the operator matrix $B=A+V$ can be block diagonalized as indicated in \eqref{eq:blockDiag1stIntro},
\eqref{eq:blockDiag2ndIntro}, and \eqref{eq:uniblockdiag} if and only if the graph subspace $\cG(\cH_0,X)$ is reducing for $A+V$
and $Y$ is a strong solution to the operator Riccati equation. A deeper discussion of various aspects of the block diagonalization
under more relaxed conditions can be found at the end of that section.

Our main result is presented in Section 3 (Theorem \ref{thm:main}). Here, we show that under mild regularity conditions the
property of the graph subspace $\cG(\cH_0,X)$ to be reducing for the operator $A+V$ and the strong solvability of the operator
Riccati equation by the skew-symmetric operator $Y$ imply one another.

In Section 4, we prove that the regularity conditions of Theorem \ref{thm:main} are automatically satisfied whenever the relative
bound of the perturbation $V$ with respect to $A$ is sufficiently small, see Theorem \ref{thm:mainRel} and Remark
\ref{rem:mainRel}. This a direct extension of the results obtained in \cite[Section 5]{AMM03}. The particular case where the
spectra of the diagonal entries $A_0$ and $A_1$ are subordinated is considered in Example \ref{examp}. The operator $A+V$ then
admits the block diagonalizations \eqref{eq:blockDiag1stIntro} and \eqref{eq:blockDiag2ndIntro} whenever the symmetric perturbation
$V$ has $A$-bound smaller than $1$. This complements the statement of the generalized Davis-Kahan $\tan2\Theta$ theorem established
in \cite{MS06} (see also \cite{GKMV13}). Finally, a more general perturbation theory for closed diagonal operators $A$ with a
suitable condition on their spectra is briefly discussed in Remark \ref{rem:moreGeneral}.

Some words about notation:

The domain of a linear operator $K$ is denoted by $\Dom(K)$ and its range by $\Ran(K)$. The restriction of $K$ to a given subset
$\cC$ of $\Dom(K)$ is written as $K|_{\cC}$.

Given another linear operator $L$, we write the extension relation $K\subset L$ (or $L\supset K$) if $L$ extends $K$, that is, if
$\Dom(K)\subset\Dom(L)$ and $Kx=Lx$ for $x\in\Dom(K)$. The operator equality $K=L$ means that $K\subset L$ and $K\supset L$.

If $K$ is a closed densely defined operator on a Hilbert space $\cH$, its adjoint operator is denoted by $K^*$ and the resolvent
set of $K$ by $\rho(K)$. The identity operator on $\cH$ is written as $I_{\cH}$. Multiples $\lambda I_{\cH}$ of the identity are
abbreviated by $\lambda$. Finally, the norm on $\cH$ is denoted by $\norm{\cdot}_{\cH}$, where the subscript $\cH$ is usually
omitted.

\vfill
\subsection*{Acknowledgements.}
The authors would like to thank Vadim Kostrykin for helpful discussions on the topic. Especially, his Ph.D.\ students
A.S.\ and S.S.\ would like to express their gratitude for introducing them to this field of research.


\section{Reducing graph subspaces and block-diagonalization}\label{sec:blockDiag}
In this section, we revisit the block diagonalization for $2\times 2$ block operator matrices as presented in
\cite[Section 5]{AMM03}. We extend the considerations to diagonally dominant symmetric operator matrices with self-adjoint diagonal
parts and provide new results towards the block diagonalization with respect to reducing graph subspaces. For a more detailed
discussion on the concept of diagonally dominant operator matrices we refer to \cite[Section 2.2]{Tre08}.

Throughout this paper we make the following assumptions:

\begin{hypothesis}\label{hyp1}
 Let $\cH=\cH_0\oplus\cH_1$ be a Hilbert space decomposed into the sum of two orthogonal subspaces $\cH_0$ and $\cH_1$. Let $A$ be
 a self-adjoint operator on $\cH$ given by the representation
 \[
  A=\begin{pmatrix} A_0 & 0\\ 0 & A_1\end{pmatrix}\,,\quad \Dom(A)=\Dom(A_0)\oplus\Dom(A_1)\,,
 \]
 with respect to the decomposition $\cH=\cH_0\oplus\cH_1$. Moreover, let $W\colon\cH_1\supset\Dom(W)\to\cH_0$ be a densely defined
 linear operator such that $\Dom(W)\supset\Dom(A_1)$ and $\Dom(W^*)\supset\Dom(A_0)$. Define the linear operator $V$ on $\cH$ by
 the $2\times 2$ block operator matrix
 \[
  V := \begin{pmatrix} 0 & W\\ W^* & 0 \end{pmatrix}\quad \text{ on }\quad \Dom(V):=\Dom(W^*) \oplus \Dom(W)\supset\Dom(A)\,.
 \]
\end{hypothesis}


\subsection{Invariant graph subspaces}

Recall that a closed subspace $\cG$ of a Hilbert space $\cH$ is said to be a \emph{graph subspace} associated with a closed
subspace $\cN\subset \cH$ and a bounded operator $X$ from $\cN$ to its orthogonal complement $\cN^\perp$ if
\[
 \cG=\cG(\cN, X):=\{x\in \cH\,|\, P_{\cN^\perp }x=XP_\cN x\}\,,
\]
where $P_\cN$ and $P_{\cN^\perp} $ denote the orthogonal projections onto $\cN$ and $\cN^\perp$, respectively.

It is easy to check that 
\[
 \cG(\cN, X)^\perp=\cG(\cN^\perp, -X^*)\,.
\]

Moreover, the orthogonal graph subspaces $\cG(\cN,X)$ and $\cG(\cN^\perp, -X^*)$ can be represented as
\begin{equation}\label{eq:graphRepr}
 \cG(\cN,X) = \Ran(T|_{\cN}) \quad\text{ and }\quad \cG(\cN^\perp, -X^*)=\Ran(T|_{\cN^\perp})\,,
\end{equation}
where the operator $T$ is given by
\[
 T = \begin{pmatrix}
      I_{\cN} & -X^*\\
      X & I_{\cN^\perp}
     \end{pmatrix}
\]
with respect to the orthogonal decomposition $\cH=\cN\oplus \cN^\perp$.

In the situation of Hypothesis \ref{hyp1}, the first step towards a block diagonalization with respect to a reducing graph subspace
is the consideration of an invariant graph subspace $\cG(\cH_0,X)$ such that its orthogonal complement is also invariant. In the
setting of unbounded operators, this requires to consider the intersections of the invariant subspaces with the operator domain. It
is therefore convenient to adopt the following notations.

\begin{hypothesis}\label{hyp2}
 Assume Hypothesis \ref{hyp1}. Let $X$ be a bounded linear operator from $\cH_0$ to $\cH_1$, and define the bounded linear
 operators $Y$ and $T$ on $\cH$ by
 \[
  Y := \begin{pmatrix} 0 & -X^*\\ X & 0 \end{pmatrix} \quad\text{ and }\quad
  T := I_{\cH}+Y=\begin{pmatrix} I_{\cH_0} & -X^*\\ X & I_{\cH_1} \end{pmatrix}\,.
 \]
 Finally, set
 \[
  \cD:=\{x\in\Dom(A) \mid Yx\in\Dom(A)\}=\cD_0\oplus\cD_1\,,
 \]
 where
 \[
  \begin{aligned}
   \cD_0 := \{f\in\Dom(A_0) \mid Xf\in\Dom(A_1)\}\,,\quad \cD_1 := \{g\in\Dom(A_1) \mid X^*g\in\Dom(A_0)\}\,.
  \end{aligned}
 \]
\end{hypothesis}
By definition, the set $\cD$ in Hypothesis \ref{hyp2} is the maximal linear subset of $\Dom (A)$ that $Y$ maps into $\Dom(A)$.

We start with the following invariance criterion, which is extracted from the proofs of \cite[Lemma 5.3 and Theorem 5.5]{AMM03}.
The corresponding reasoning is taken over and repeated briefly.

\begin{lemma}[cf.\ {\cite[Lemma 5.3 and Theorem 5.5]{AMM03}}]\label{lem:invariance}
 Assume Hypotheses \ref{hyp1} and \ref{hyp2}.
 
 The following are equivalent:
 \begin{enumerate}
  \renewcommand{\theenumi}{\roman{enumi}}
  \item The graph subspaces $\cG(\cH_0,X)$ and $\cG(\cH_1,-X^*)$ are invariant for the operator $A+V$.
  \item The operator $Y$ satisfies the Riccati equation
        \begin{equation}\label{eq:riccati}
         AYx-YAx-YVYx+Vx=0 \quad\text{ for }\quad x\in\cD\,.
        \end{equation}
  \item The operator $T$ satisfies
        \begin{equation}
         (A+V)Tx = T(A+VY)x\quad\text{ for }\quad x\in\cD\,.\label{eq:qsimilar(A+V)T}
        \end{equation}
  \item The operator $T^*$ satisfies
        \begin{equation}\label{eq:qsimilarT*(A+V)}
         T^*(A+V)x = (A-YV)T^*x\,\quad\text{ for }\quad x\in\cD\,.
        \end{equation}
 \end{enumerate}

 \begin{proof} It is easy to see that the graph subspace $\cG(\cH_0,X)$ is invariant for the operator $A+V$ if and only if the
  operator $X$ satisfies the equation
  \begin{equation}
   A_1Xf - XA_0f - XWXf + W^*f = 0 \quad\text{ for }\quad f\in\cD_0\,.\label{eq:riccatiSplitX}
  \end{equation}
  
  Indeed, observing that
  \[
   \Dom(A)\cap\cG(\cH_0,X) = \{f\oplus Xf \mid f\in\cD_0\}
  \]
  and that
  \[
   (A+V)\begin{pmatrix} f\\ Xf \end{pmatrix} = \begin{pmatrix} A_0f + WXf\\ W^*f + A_1Xf \end{pmatrix} \quad\text{ for }\quad
   f\in\cD_0\,,
  \]
  it follows that the graph subspace $\cG(\cH_0,X)$ is invariant for $A+V$ if and only if the equation $(W^* + A_1X)f=X(A_0 + WX)f$
  holds for all $f\in\cD_0$, which, in turn, can be rewritten as \eqref{eq:riccatiSplitX}.
  
  Upon observing that
  \[
   \Dom(A)\cap\cG(\cH_1,-X^*) = \{-X^*g\oplus g \mid g\in\cD_1\}\,,
  \]
  one concludes in a completely analogous way that the graph subspace $\cG(\cH_1,-X^*)$ is invariant for the operator $A+V$ if and
  only if $X^*$ satisfies the equation
  \begin{equation}
   A_0X^*g - X^*A_1g + X^*W^*X^*g - Wg = 0 \quad\text{ for }\quad g\in\cD_1\,.\label{eq:riccatiSplitX*}
  \end{equation}
  
  Thus, both $\cG(\cH_0,X)$ and its orthogonal complement $\cG(\cH_0,X)^\perp=\cG(\cH_1,-X^*)$ are invariant for the operator $A+V$
  if and only if the pair of Riccati equations \eqref{eq:riccatiSplitX} and \eqref{eq:riccatiSplitX*} hold. Now, it is easy to
  check that this pair of equations can be rewritten as the single Riccati equation \eqref{eq:riccati}. This proves the stated
  equivalence of (i) and (ii).
  
  Finally, taking into account that $T=I_{\cH}+Y$ and $T^*=I_{\cH}-Y$, the two equations \eqref{eq:qsimilar(A+V)T} and
  \eqref{eq:qsimilarT*(A+V)} are just reformulations of the Riccati equation \eqref{eq:riccati}, so that (iii) and (iv) are
  equivalent to (ii).
 \end{proof}%
\end{lemma}

\begin{remark}\label{rem:invariance}
 In essence, equation \eqref{eq:qsimilar(A+V)T} in Lemma \ref{lem:invariance} has already been considered in
 \cite[Theorem 5.5 (ii)]{AMM03}, while \eqref{eq:qsimilarT*(A+V)} has not been mentioned in \cite{AMM03}.
 
 However, both equations seem to have their own right. The one \eqref{eq:qsimilar(A+V)T} appears naturally from the representation
 \eqref{eq:graphRepr}, but \eqref{eq:qsimilarT*(A+V)} proves to be more useful for some of our considerations, in particular those
 in Section \ref{sec:mainResult}. One reason for this is that the operator $A-YV$ has natural domain $\Dom(A)$, whereas the natural
 domain of the operator $A+VY$ satisfies
 \[
  \cD \subset \Dom(A+VY) \subset \Dom(A)\,,
 \]
 and either one of the inclusions may a priori be strict. In particular, $\Dom(A+VY)$ is determined by the choice of $Y$. We
 therefore focus on equation \eqref{eq:qsimilarT*(A+V)} in the following considerations. This is discussed in more detail at the
 end of this section.
\end{remark}

Before we discuss reducing graph subspaces for the operator $A+V$, a more detailed study of the operators $T$ and $T^*$ is
required.


\subsection{The operators $T$ and $T^*$}
We start with the following elementary observation that relates some mapping properties of $T$, $T^*$, and $Y$.

\begin{lemma}\label{lem:TT*Y}
 Assume Hypotheses \ref{hyp1} and \ref{hyp2}, and let $\cC$ be a linear subset of $\Dom(A)$. The following are equivalent:
 \begin{enumerate}
  \renewcommand{\theenumi}{\roman{enumi}}
  \item $T$ maps $\cC$ into $\Dom(A)$.
  \item $T^*$ maps $\cC$ into $\Dom(A)$.
  \item $Y$ maps $\cC$ into $\Dom(A)$.
 \end{enumerate}
 In particular, $\cD$ is the maximal linear subset of $\Dom(A)$ that $T$ (resp.\ $T^*$) maps into $\Dom(A)$.
 \begin{proof}
  Taking into account that $T=I_{\cH}+Y$ and $T^*=I_{\cH}-Y$, the equivalence is obvious. The additional statement follows from the
  fact that $\cD$ is, by definition, the maximal linear subset of $\Dom(A)$ that $Y$ maps into $\Dom(A)$.
 \end{proof}%
\end{lemma}

It is easy to see that the operators $T$ and $T^*$ each have a bounded inverse. Indeed, the spectrum of $Y$ is a subset of the
imaginary axis since $Y^*=-Y$. Hence, zero belongs to the resolvent sets of $T=I_{\cH}+Y$ and $T^*=I_{\cH}-Y$
(cf.\ \cite[Theorem 5.5 (i)]{AMM03}).

It turns out that the mapping properties of $T^{-1}$ and $(T^*)^{-1}$ are related in a similar way as the ones of $T$ and $T^*$.
\begin{lemma}\label{lem:TT*inv}
 Assume Hypotheses \ref{hyp1} and \ref{hyp2}. Then, $T^{-1}$ maps $\Dom(A)$ into $\Dom(A)$ if and only if $(T^*)^{-1}$ does.
 
 In this case, the set $\cD$ can alternatively be represented as
 \begin{equation}\label{eq:altReprD}
  \cD = \Ran(T^{-1}|_{\Dom(A)}) = \Ran\bigl((T^*)^{-1}|_{\Dom(A)}\bigr)\,.
 \end{equation}

 \begin{proof}
  Introduce the unitary block diagonal matrix
  \begin{equation}\label{eq:defJ}
   J = \begin{pmatrix}
        I_{\cH_0} & 0\\
        0 & -I_{\cH_1}
       \end{pmatrix}
  \end{equation}
  with respect to the orthogonal decomposition $\cH=\cH_0\oplus\cH_1$. Obviously, $J=J^*$ maps $\Dom(A)$ onto itself, and one has
  $T^*=JTJ^*$, so that
  \[
   (T^*)^{-1}=JT^{-1}J^*\,.
  \]
  Hence, $T^{-1}$ maps $\Dom(A)$ into $\Dom(A)$ if and only if $(T^*)^{-1}$ does.

  In this case, $T$ maps $\Ran\bigl(T^{-1}|_{\Dom(A)}\bigr)\subset\Dom(A)$ into $\Dom(A)$. Since $\cD$ is by Lemma \ref{lem:TT*Y}
  the maximal linear subset of $\Dom(A)$ that $T$ maps into $\Dom(A)$, one concludes that
  $\Ran\bigl(T^{-1}|_{\Dom(A)}\bigr)\subset\cD$, that is,
  \[
   \Dom(A)\subset\Ran(T|_\cD)\,.
  \]
  Since also $\Ran(T|_\cD)\subset\Dom(A)$ by definition of $\cD$, this yields $\Ran(T|_{\cD})=\Dom(A)$, which is equivalent to the
  first equality in \eqref{eq:altReprD}. The analogous reasoning for $(T^*)^{-1}$ shows that also
  \[
   \Dom(A)=\Ran(T^*|_\cD)\,,
  \]
  which proves the second equality in \eqref{eq:altReprD}.
 \end{proof}%
\end{lemma}


\subsection{Reducing graph subspaces}  
Our next result shows that the equivalence in Lemma \ref{lem:TT*inv} automatically takes place if the graph subspace $\cG(\cH_0,X)$
is reducing for $A+V$. In fact, we have the following criterion for the graph subspace $\cG(\cH_0,X)$ to be reducing for the
operator $A+V$.

\begin{lemma}\label{lem:reducing}
 Assume Hypotheses \ref{hyp1} and \ref{hyp2}. The following are equivalent:
 \begin{enumerate}
  \renewcommand{\theenumi}{\roman{enumi}}
  \item The graph subspace $\cG(\cH_0,X)$ is reducing for the operator $A+V$, that is, $\cG(\cH_0,X)$ and $\cG(\cH_1,-X^*)$ are
        invariant for $A+V$ and $\Dom(A+V)=\Dom(A)$ splits as
        $\Dom(A)=\bigl(\Dom(A)\cap\cG(\cH_0,X)\bigr)+\bigl(\Dom(A)\cap\cG(\cH_1,-X^*)\bigr)$.
  \item The graph subspaces $\cG(\cH_0,X)$ and $\cG(\cH_1,-X^*)$ are invariant for $A+V$, and $T^{-1}$ (resp.\ $(T^*)^{-1}$) maps
        $\Dom(A)$ into itself.
  \item One has
        \begin{equation}\label{eq:reducingLeftExt}
         T^*(A+V) \supset (A-YV)T^*\,.
        \end{equation}
 \end{enumerate}

 \begin{proof} Clearly, one has
  \begin{equation*}
  \Dom(A+V)\cap\cG(\cH_0,X) = \{f\oplus Xf \mid f\in\cD_0\}=\Ran(T|_{\cD_0})\,,
  \end{equation*}
  as well as
  \begin{equation*}
   \Dom(A+V)\cap\cG(\cH_1,-X^*) = \{-X^*g\oplus g \mid g\in\cD_1\}=\Ran(T|_{\cD_1})\,.
  \end{equation*}
  This yields
  \begin{equation*}
   \bigl(\Dom(A+V)\cap\cG(\cH_0,X)\bigr) + \bigl(\Dom(A+V)\cap\cG(\cH_1,-X^*)\bigr) = \Ran(T|_{\cD}) \subset\Dom(A)\,.
  \end{equation*}

  Therefore, the graph subspace $\cG(\cH_0,X)$ is reducing for $A+V$ if and only if $\cG(\cH_0,X)$ and its orthogonal complement
  $\cG(\cH_1,-X^*)$ are invariant for $A+V$ and
  \begin{equation}\label{eq:reducingT}
   \Dom(A) = \Dom(A+V) = \Ran(T|_\cD)\,.
  \end{equation}
  It follows from Lemma \ref{lem:TT*inv} that \eqref{eq:reducingT} holds if and only if $T^{-1}$ (resp.\ $(T^*)^{-1}$) maps
  $\Dom(A)$ into $\Dom(A)$. This proves the equivalence of (i) and (ii).

  Suppose that (i) holds. Then, by Lemma \ref{lem:invariance} the invariance of the graph subspaces $\cG(\cH_0,X)$ and
  $\cG(\cH_1,-X^*)$ implies that
  \[
   T^*(A+V)x = (A-YV)T^*x \quad\text{ for }\quad x\in\cD\,.
  \]
  Moreover, together with \eqref{eq:reducingT} one also has
  \begin{equation}\label{eq:reducingT*}
   \Dom(A) = \Ran(T^*|_\cD)\,,
  \end{equation}
  which is due to Lemma \ref{lem:TT*inv}. Now, it follows from \eqref{eq:reducingT*} that
  \[
   \Dom\bigl((A-YV)T^*\bigr)=\cD\subset\Dom(A+V)=\Dom\bigl(T^*(A+V)\bigr)
  \]
  since $\Dom(A-YV)=\Dom(A)$. Hence, (i) implies the extension relation (iii).

  Conversely, suppose that (iii) holds. Since $\cD\subset\Dom\bigl((A-YV)T^* \bigr)$, this implies that
  \[
   T^*(A+V)x = (A-YV)T^*x\quad\text{ for }\quad x\in\cD\,.
  \]
  Hence, the graph subspaces $\cG(\cH_0,X)$ and $\cG(\cH_1,-X^*)$ are invariant for $A+V$ by Lemma \ref{lem:invariance}. Moreover,
  \[
   \Ran\bigl((T^*)^{-1}|_{\Dom(A)}\bigr) = \Dom((A-YV)T^*) \subset \Dom\bigl(T^*(A+V)\bigr) = \Dom(A)\,,
  \]
  so that (ii) holds.

  This completes the proof.
 \end{proof}%
\end{lemma}

\begin{remark}\label{rem:sandwich}
 In the situation of Lemma \ref{lem:reducing}, an analogous reasoning shows that the graph subspace $\cG(\cH_0,X)$ is reducing for
 the operator $A+V$ if and only if
 \[
  (A+V)T \subset T(A+VY)\,.
 \]
 Hence, the property of $\cG(\cH_0,X)$ to be reducing for $A+V$ can be characterized by the following chain of operator extension
 relations:
 \begin{equation}\label{eq:sandwich}
  (T^*)^{-1} (A-YV)T^*\subset   (A+V) \subset T(A+VY)T^{-1}\,.
 \end{equation}
 For a more detailed discussion of this ``sandwich'' property \eqref{eq:sandwich} we refer to Subsection \ref{subsec:diss} at the
 end of this section. For now, it suffices to observe that \eqref{eq:sandwich} becomes a series of operator equalities if
 $\cD=\Dom(A)$. Indeed, in this case one has $\Dom(A+VY)=\Dom(A)$, and $T$ and $T^*$ both map $\Dom(A)$ onto $\Dom(A)$ by Lemmas
 \ref{lem:TT*Y} and \ref{lem:TT*inv}. In fact, the converse statement also holds, see Theorem \ref{thm:blockDiag}.
\end{remark}


\subsection{Strong solutions to the operator Riccati equation}
Starting from the extension relation \eqref{eq:reducingLeftExt} in Lemma \ref{lem:reducing}, the discussion in Remark
\ref{rem:sandwich} shows that more information on the set $\cD$ is needed in order to get the corresponding operator equality in
\eqref{eq:reducingLeftExt}. The following result states that such kind of information is encoded in the strong solvability of the
operator Riccati equation by the operator $Y$.

\begin{lemma}\label{lem:riccati}
 Assume Hypotheses \ref{hyp1} and \ref{hyp2}. Then, the extension relation
 \[
  T^*(A+V) \subset (A-YV)T^*
 \]
 holds if and only if $Y$ is a strong solution to the operator Riccati equation
 \[
  AY - YA - YVY + V = 0\,,
 \]
 that is, $\Ran(Y|_{\Dom(A)})\subset\Dom(A)$ and
 \[
  AYx - YAx - YVYx + Vx = 0 \quad\text{ for }\quad x\in\Dom(A)\,.
 \]
 In this case, one has $\cD=\Dom(A)$.         
 
 \begin{proof}
  By Lemma \ref{lem:TT*Y}, the operator $Y$ maps $\Dom(A)$ into $\Dom(A)$ if and only if $T^*$ does. Moreover, one has
  $\Dom(A)=\cD$ in this case since $\cD$ is the maximal linear subset of $\Dom(A)$ that $Y$ maps into $\Dom(A)$. Taking into
  account that
  \[
   \Dom(A+V)=\Dom(A)=\Dom(A-YV)\,,
  \]
  the stated equivalence now follows from Lemma \ref{lem:invariance}.
 \end{proof}%
\end{lemma}


\subsection{Block diagonalization}
We are now ready to present the main result of this section. It yields necessary and sufficient conditions under which equation
\eqref{eq:qsimilarT*(A+V)} in Lemma \ref{lem:invariance} extends to an operator equality. As the theorem states, in this case also
\eqref{eq:qsimilar(A+V)T} becomes an operator equality.
\begin{theorem}\label{thm:blockDiag}
 Assume Hypotheses \ref{hyp1} and \ref{hyp2}. Then, the operator $A+V$ admits the block diagonalization
 \begin{equation}\label{eq:blockDiag2nd}
  T^*(A+V) (T^*)^{-1}  = A-YV = \begin{pmatrix} A_0 + X^*W^* & 0\\ 0 & A_1-XW\end{pmatrix}
 \end{equation}
 if and only if
 \begin{enumerate}
  \renewcommand{\theenumi}{\roman{enumi}}
  \item the graph subspace $\cG(\cH_0,X)$ is reducing for $A+V$,

  and

  \item the operator $Y$ is a strong solution to the
        operator Riccati equation
        \[
         AY - YA - YVY + V = 0\,.
        \]
 \end{enumerate}
 In this case, $A+V$ also admits the block diagonalization 
 \begin{equation}\label{eq:blockDiag1st}
  T^{-1}(A+V) T  = A+VY = \begin{pmatrix} A_0 + WX & 0\\ 0 & A_1-W^*X^*\end{pmatrix}\,.
 \end{equation}
 
 \begin{proof}
  The stated equivalence is just a combination of Lemmas \ref{lem:reducing} and \ref{lem:riccati}.
   
  In this case, one has $\cD=\Dom(A)$, and $T$ and $T^*$ both map $\Dom(A)$ onto $\Dom(A)$ by Lemmas \ref{lem:TT*Y} and
  \ref{lem:TT*inv}. In particular, it follows that
  \[
   \cD = \Dom(A+VY) = \Dom(A)\,.
  \]
  The block diagonalization \eqref{eq:blockDiag1st} is then a direct consequence of equation \eqref{eq:qsimilar(A+V)T} in Lemma
  \ref{lem:invariance}.
 \end{proof}%
\end{theorem}

\begin{remark}\label{rem:unitary}
 In the situation of Theorem \ref{thm:blockDiag}, one has
 \begin{equation}\label{eq:T*T}
  T^*T  (A+VY) (T^*T)^{-1} = A-Y V
 \end{equation}
 by \eqref{eq:blockDiag2nd} and \eqref{eq:blockDiag1st}. In particular, the operators $A+V$, $A+VY$, and $A-YV$ are similar to one
 another. Moreover, the operator $T$ is normal and $T^*T$ is block diagonal,
 \[
  T^*T= \begin{pmatrix} I+X^*X & 0\\0& I+XX^* \end{pmatrix}= TT^*\,.
 \]
 Hence, the corresponding entries of the block diagonal matrices $A+VY$ and $A-YV$ are also similar to one another. More
 explicitly, 
 \[
  (I+X^*X)(A_0 + WX) (I+X^*X)^{-1}=A_0+X^*W^*
 \]
 and
 \[
  (I+XX^*)(A_1-W^*X^*)(I+XX^*)^{-1}=A_1-XW\,.
 \]
 In addition, the operator $A+V$ is unitarily equivalent to a block diagonal matrix via the unitary transformation $U$ from the 
 polar decomposition $T=U\abs{T}$ (cf.\ \cite[Theorem 5.5 (iii)]{AMM03}). More precisely, together with
 $T^*=\abs{T} U^*=\abs{T^*} U^*$, one has
 \begin{equation}\label{eq:unitaryBlockDiag}
  U^*(A+V)U = \abs{T}(A+VY)\abs{T}^{-1} = \abs{T}^{-1}(A-YV)\abs{T} = \begin{pmatrix} B_0 & 0\\ 0 & B_1 \end{pmatrix}\,,
 \end{equation}
 with
 \[
  B_0:=(I+X^*X)^{1/2}(A_0+WX)(I+X^*X)^{-1/2}
 \]
 and
 \[
  B_1:=(I+XX^*)^{1/2}(A_1-W^*X^*)(I+XX^*)^{-1/2}
 \]
 on the corresponding natural domains
 \[
  \Dom(B_0):=\Ran\bigl((I+X^*X)^{1/2}|_{\Dom(A_0)}\bigr)
 \]
 and
 \[
  \Dom(B_1):=\Ran\bigl((I+XX^*)^{1/2}|_{\Dom(A_1)}\bigr)\,.
 \]

 Note that the second equality in \eqref{eq:unitaryBlockDiag} is due to \eqref{eq:T*T}. In particular, the two block
 diagonalizations \eqref{eq:blockDiag2nd} and \eqref{eq:blockDiag1st} lead to the same unitary block diagonalization
 \eqref{eq:unitaryBlockDiag}.
\end{remark}

The property of the operators $T$ and $T^*$ to map $\Dom(A)$ onto $\Dom(A)$ is essential for the block diagonalizations
\eqref{eq:blockDiag2nd} and \eqref{eq:blockDiag1st} to hold. Moreover, by Lemmas \ref{lem:TT*Y} and \ref{lem:TT*inv}, the operator
$T$ has this property if and only if $T^*$ does. Following the proof of Lemma \ref{lem:TT*inv}, a crucial point here is the fact
that the two operators $T$ and $T^*=JTJ^*$ are coupled by the operator $J$ defined by \eqref{eq:defJ}. This operator $J$ also
satisfies
\[
 J(A+V)J^*=A+JVJ^*=A-V\,.
\]
In particular, the operators $A\pm V$ can be diagonalized simultaneously, and the operator $T^*$ is related to $A-V$ in the same
way as $T$ is to $A+V$. This observation generalizes to the following statement.

\begin{remark}
 For $\theta$ from the complex unit circle introduce the unitary block diagonal matrix
 \[
  J_\theta := \begin{pmatrix}
               I_{\cH_0} & 0\\ 0 & \theta\cdot I_{\cH_1}
              \end{pmatrix}\,.
 \]
 Then, $J_\theta$ and $J_\theta^*=J_{\overline{\theta}}$ map $\Dom(A)$ onto itself, and one has
 \[
  J_\theta (A+V)J_\theta^* = A + J_\theta VJ_\theta^* =
  \begin{pmatrix}
   A_0 & \overline{\theta}\cdot W\\
   \theta\cdot W^* & A_1
  \end{pmatrix}\,.
 \]
 The invariant and reducing graph subspaces for $A+V$ transform accordingly to invariant and reducing graph subspaces for
 $A+J_{\theta}VJ_{\theta}^*$, respectively. In this case, $X$ and $W$ are replaced by $\theta\cdot X$ and
 $\overline{\theta}\cdot W$, respectively. Hence, all the operators $A+J_\theta VJ_\theta^*$, $\abs{\theta}=1$, can be diagonalized
 simultaneously, where the corresponding diagonal operators are the same for all $\theta$ since diagonal operators are invariant
 under conjugation by $J_\theta$. More precisely,
 \[
  A + (J_\theta VJ_\theta^*)(J_\theta YJ_\theta^*) = A+VY \quad\text{ and }\quad
  A - (J_\theta YJ_\theta^*)(J_\theta VJ_\theta^*) = A-YV\,.
 \]
\end{remark}


\subsection{Further discussion and more general statements}\label{subsec:diss}

We now give an informal but more detailed discussion on reducing graph subspaces and the related block diagonalizations. We also
briefly discuss some more general situations to which our considerations can be carried over. Note that the content of this
subsection is not required for Sections \ref{sec:mainResult} and \ref{sec:relBddPert}. However, some remarks in these sections pick
up the following discussions.

Recall that by Lemma \ref{lem:reducing} and the discussion in Remark \ref{rem:sandwich} a reducing graph subspace $\cG(\cH_0,X)$
for the operator $A+V$ can be characterized by the chain of operator extension relations
\begin{equation}\label{eq:sandwichDiss}
 (T^*)^{-1}(A-YV)T^* \subset A+V \subset T(A+VY)T^{-1}\,.
\end{equation}
Moreover, as described in Remark \ref{rem:invariance}, one has
\[
 \Dom(A+V) = \Dom(A) = \Dom(A-YV)\,,
\]
whereas the natural domain of the operator $A+VY$ satisfies
\begin{equation}\label{eq:DIncl}
 \cD \subset \Dom(A+VY) \subset \Dom(A)\,.
\end{equation}
Note that either one of the extensions in \eqref{eq:sandwichDiss}, as well as of the inclusions in \eqref{eq:DIncl}, may a priori
be strict.

Taking into account the representation of the set $\cD$ from Lemma \ref{lem:TT*inv},
\begin{equation}\label{eq:reprD}
 \cD = \Ran(T^{-1}|_{\Dom(A)}) = \Ran\bigl((T^*)^{-1}|_{\Dom(A)}\bigr)\,,
\end{equation}
the right-hand extension relation in \eqref{eq:sandwichDiss} yields a block diagonalization for the operator $A+V$,
\begin{equation}\label{eq:blockDiagTD}
 T^{-1}(A+V)T = (A+VY)|_{\cD}\,.
\end{equation}
The left-hand relation in \eqref{eq:sandwichDiss}, however, allows for a block diagonalization only of the restricted operator
$(A+V)|_{\cD}$, that is,
\begin{equation}\label{eq:blockDiagDT*}
 T^*(A+V)|_{\cD}(T^*)^{-1} = A-YV\,.
\end{equation}

Now, the right-hand relation in \eqref{eq:sandwichDiss} becomes an equality if and only if $T$ maps the domain $\Dom(A+VY)$ onto
$\Dom(A)$. Taking into account \eqref{eq:DIncl}, \eqref{eq:reprD}, and Lemma \ref{lem:TT*Y}, this, in turn, is equivalent to
$\cD=\Dom(A+VY)$. In this case, we say that $A+V$ admits \emph{the first block diagonalization}
\begin{equation}\label{eq:blockDiag1stDiss}
 T^{-1}(A+V)T = A+VY\,.
\end{equation}
If, in addition, $\Dom(A+VY) \subsetneq \Dom(A)$, then $\cD\subsetneq\Dom(A)$ as well, and $Y$ is not a strong solution to the
operator Riccati equation $AY-YA-YVY+V=0$. In particular, the left-hand relation in \eqref{eq:sandwichDiss} is then strict by
Theorem \ref{thm:blockDiag}.

The left-hand relation in \eqref{eq:sandwichDiss} becomes an equality if and only if $T^*$ maps $\Dom(A)$ onto $\Dom(A)$. Since
$(T^*)^{-1}$ by Lemma \ref{lem:TT*inv} already maps $\Dom(A)$ into $\Dom(A)$, it follows from Lemma \ref{lem:TT*Y} that this is
equivalent to $\cD=\Dom(A)$. In this case, we say that $A+V$ admits \emph{the second block diagonalization}
\begin{equation}\label{eq:blockDiag2ndDiss}
 T^*(A+V)(T^*)^{-1} = A-YV\,.
\end{equation}
By \eqref{eq:DIncl} one also has $\cD=\Dom(A+VY)=\Dom(A)$, so that $A+V$ admits the first block diagonalization
\eqref{eq:blockDiag1stDiss} as well.

In summary, for a reducing graph subspace $\cG(\cH_0,X)$ for $A+V$ there are three cases that need to be distinguished:
\begin{enumerate}
 \renewcommand{\theenumi}{\alph{enumi}}
 \item $\cD\subsetneq\Dom(A+VY)$. In this case, $A+V$ admits only the block diagonalization \eqref{eq:blockDiagTD} with the
       diagonal operator $A+VY$ being restricted to $\cD$. Both extension relations in \eqref{eq:sandwichDiss} are strict.
 \item $\cD=\Dom(A+VY)\subsetneq\Dom(A)$. In this case, $A+V$ admits the first block diagonalization \eqref{eq:blockDiag1stDiss},
       but not the second \eqref{eq:blockDiag2ndDiss}. The restriction $(A+V)|_{\cD}$ satisfies \eqref{eq:blockDiagDT*}.
 \item $\cD=\Dom(A)$. In this case, $A+V$ admits the second block diagonalization \eqref{eq:blockDiag2ndDiss} and, hence, also the
       first \eqref{eq:blockDiag1stDiss}. Moreover, the operator $Y$ is a strong solution to the operator Riccati equation.
\end{enumerate}

The linear subsets $\cD$ and $\Dom(A+VY)$ depend on the choice of the operator $Y$ and additional information on them may not be
available in advance. Both sets are therefore given only on an abstract level, which makes the cases (a) and (b) more difficult to
handle in applications. Moreover, in this general situation we neither are able to give an example that the cases (a) and (b)
really do occur nor can we prove otherwise. For these reasons we have focused on case (c) throughout this section. For a particular
situation where only the case (c) occurs, we refer to Section \ref{sec:relBddPert}.

However, case (b) deserves a separate discussion. More precisely, we obtain the following formal criterion for the first block
diagonalization \eqref{eq:blockDiag1stDiss} to hold. This result is analogous to Theorem \ref{thm:blockDiag}. With the above
considerations, the corresponding proof is straightforward and is hence omitted.
\begin{theorem}\label{thm:blockDiag1st}
 Assume Hypotheses \ref{hyp1} and \ref{hyp2}. Then, the operator $A+V$ admits the first block diagonalization
 \[
  T^{-1}(A+V) T  = A+VY = \begin{pmatrix} A_0 + WX & 0\\ 0 & A_1-W^*X^*\end{pmatrix}
 \]
 if and only if
 \begin{enumerate}
  \renewcommand{\theenumi}{\roman{enumi}}
  \item the graph subspace $\cG(\cH_0,X)$ is reducing for $A+V$,

  and

  \item the operator $Y$ satisfies $\Ran(Y|_{\Dom(A+VY)})\subset\Dom(A)$ and
        \[
         AYx - YAx - YVYx + Vx = 0 \quad\text{ for }\quad x\in\Dom(A+VY)\,.
        \]
 \end{enumerate}
\end{theorem}

We close this section with two concluding remarks that address the extension of the results in this section to certain more general
situations.

The first one considers the case where the operator $A$ is not assumed to be self-adjoint. 
\begin{remark}\label{rem:selfAdjoint}
 The hypothesis that $A$ is a self-adjoint operator has only a formal character and is nowhere used explicitly. Thus, all results
 in this section, especially Theorem \ref{thm:blockDiag}, remain valid if $A$ is assumed to be just a linear operator with
 $\Dom(A)\subset \Dom(V)$. However, if $A$ is not self-adjoint (or not even symmetric), the a priori consideration of orthogonal
 invariant graph subspaces for the operator $A+V$ seems to be unmotivated for those may not exist in principle. The assumption that
 they do exist may distract from the main issue of this work. That is why only the case of self-adjoint operators $A$ has been
 discussed.
\end{remark}

The second remark discusses perturbations $V$ for which the hypothesis $\Dom(A)\subset\Dom(V)$ is not satisfied.

\begin{remark}\label{rem:general}
 The hypothesis $\Dom(A)\subset \Dom(V)$ can also be dropped throughout this section. In a more general sense, $Y$ is said to be a
 strong solution to the operator Riccati equation $AY-YA-YVY+V=0$ if
 \[
  \Ran(Y|_{\Dom(A)\cap\Dom(V)})\subset \Dom(A)\cap\Dom(V)\,,
 \]
 and
 \[
  AYx-YAx-YVYx+Vx=0 \quad\text{ for }\quad x\in\Dom(A)\cap\Dom(V)\,.
 \]
 Since
 \[
  \Dom(A+V) = \Dom(A)\cap\Dom(V) = \Dom(A-YV)\,,
 \]
 the domain $\Dom(A)$ then has to be replaced by $\Dom(A)\cap\Dom(V)$ everywhere. For example, the set $\cD$ from Hypothesis
 \ref{hyp2} has to be defined as
 \[
  \cD:=\{ x\in\Dom(A)\cap\Dom(V) \mid Yx\in\Dom(A)\cap\Dom(V)\}\,.
 \]
 This allows, for instance, to consider unbounded perturbations $V$ if the diagonal operator $A$ is bounded. However, if
 $\Dom(A)\not\subset\Dom(V)$, the technique in Section \ref{sec:relBddPert} can not be applied directly, so that the assumption
 $\Dom(A)\subset\Dom(V)$ seems to be reasonable for our purposes.
\end{remark}


\section{Reducing graph subspaces and strong solutions to Riccati equations}\label{sec:mainResult}

In Theorem \ref{thm:blockDiag}, the two block diagonlizations \eqref{eq:blockDiag2nd} and \eqref{eq:blockDiag1st} for $A+V$ hold if
and only if $\cG(\cH_0,X)$ is a reducing subspace for $A+V$ and $Y$ is a strong solution to the operator Riccati equation. It is a
natural question under which (additional) assumptions one of these conditions implies the other and is therefore sufficient for the
block diagonalization \eqref{eq:blockDiag2nd} (and hence also \eqref{eq:blockDiag1st}) to hold. This matter is investigated in the
present section.

We start with the following elementary observation.
\begin{lemma}[{\cite[Lemma 1.3]{Schm12}}]\label{lem:schm}
 Let $\cT$ and $\cS$ be linear operators such that $\cS\subset\cT$. If $\cS$ is surjective and $\cT$ is injective, then $\cS=\cT$.
 \begin{proof}
  For the sake of completeness, we reproduce the proof from \cite{Schm12}.
 
  Let $y\in\Dom(\cT)$ be arbitrary. Since $\cS$ is surjective, there is $x\in\Dom(\cS) \subset\Dom(\cT)$ such that
  $\cT y=\cS x=\cT x$, where we have taken into account that $\cS\subset\cT$. The injectivity of $\cT$ now implies that
  $y=x\in\Dom(\cS)$. Thus, $\Dom(\cT)=\Dom(\cS)$ and, hence, $\cS=\cT$.
 \end{proof}%
\end{lemma}
In view of the extension relations in Lemmas \ref{lem:reducing} and \ref{lem:riccati}, we use the preceding lemma in the following
form.
\begin{corollary}\label{cor:schm}
 Let $K$ and $L$ be linear operators on Hilbert spaces $\cK$ and $\cL$, respectively, and let $S\colon\cK\to\cL$ be an isomorphism.
 Suppose that
 \[
  SK\subset LS\,.
 \]
 
 If there exists some constant $\lambda$ such that $K-\lambda$ is surjective and $L-\lambda$ is injective, then
 \[
  SK=LS
 \]
 holds as an operator equality. In this case, both operators $K-\lambda$ and $L-\lambda$ are bijective.
 
 \begin{proof}
  By a simple shift argument, it suffices to assume that $\lambda=0$, that is, $K$ is surjective and $L$ is injective. Since $S$ is
  an isomorphism, the operator $SK$ also is surjective and $LS$ is injective.

  The claim that $SK=LS$ holds as an operator equality now follows from Lemma \ref{lem:schm}. The additional statement is then
  obvious.
 \end{proof}%
\end{corollary}

\begin{remark}\label{rem:regCondClosed}
 If, in addition to the hypothesis of Corollary \ref{cor:schm}, the operators $K$ and $L$ are assumed to be closed and their
 resolvent sets are not disjoint, that is, $\rho(K)\cap\rho(L)\neq\emptyset$, then the intertwining relation $SK\subset LS$
 automatically implies the equality $SK=LS$. Indeed, for every $\lambda\in\rho(K)\cap\rho(L)$, the operators $K-\lambda$ and
 $L-\lambda$ are bijective.
\end{remark}

The following theorem constitutes the central result of this work. It states that either one of the conditions (i) and (ii) in
Theorem \ref{thm:blockDiag} implies the other if certain regularity conditions on $A+V$ and $A-YV$ as operators from $\Dom(A)$ to
$\cH$ are satisfied.

\begin{theorem}\label{thm:main}
 Assume Hypotheses \ref{hyp1} and \ref{hyp2}.
 \begin{enumerate}
  \renewcommand{\theenumi}{\alph{enumi}}
  \item Let the graph subspace $\cG(\cH_0,X)$ be reducing for the operator $A+V$. Then, $Y$ is a strong solution to the
        operator Riccati equation
        \[
         AY-YA-YVY+V=0
        \]
        if there is some constant $\lambda$ such that $A+V-\lambda$ is injective and $A-YV-\lambda$ is surjective.
  \item Let $Y$ be a strong solution to the operator Riccati equation $AY-YA-YVY+V=0$. Then, the graph subspace
        $\cG(\cH_0,X)$ is reducing for $A+V$, if there is some constant $\lambda$ such that $A+V-\lambda$ is surjective and
        $A-YV-\lambda$ is injective.
 \end{enumerate}
 \begin{proof}
  (a). Since $\cG(\cH_0,X)$ is reducing for $A+V$ by hypothesis, it follows from Lemma \ref{lem:reducing} that
  \[
   (A-YV)T^* \subset T^*(A+V)\,,
  \]
  which can be rewritten as
  \[
   (T^*)^{-1}(A-YV) \subset (A+V)(T^*)^{-1}\,.
  \]
  Under the given assumptions, it follows from Corollary \ref{cor:schm} that the identity
  \[
   (T^*)^{-1}(A-YV) = (A+V)(T^*)^{-1}
  \]
  holds as an operator equality, which, in turn, can be rewritten as the block diagonalization \eqref{eq:blockDiag2nd}. This proves
  the claim by Theorem \ref{thm:blockDiag}.
  
  (b). Since $Y$ is a strong solution to the operator Riccati equation, it follows from  Lemma \ref{lem:riccati} that
  \[
   T^*(A+V) \subset (A-YV)T^*\,.
  \]
  Now, Corollary \ref{cor:schm} implies that the identity
  \[
   T^*(A+V) = (A-YV)T^*
  \]
  holds as an operator equality, which can be rewritten as \eqref{eq:blockDiag2nd}. Thus, the claim follows again from Theorem
  \ref{thm:blockDiag}.
 \end{proof}%
\end{theorem}

\begin{remark}
 It is interesting to note that the regularity conditions in statement (b) of Theorem \ref{thm:main} are switched compared to those
 in statement (a). More precisely, in (a) it is assumed that the operator $A+V-\lambda$ is injective and $A-YV-\lambda$ is
 surjective while in statement (b) the operator $A+V-\lambda$ is supposed to be surjective and $A-YV-\lambda$ to be injective.
 However, due to Corollary \ref{cor:schm}, in both parts of Theorem \ref{thm:main} the operators $A+V-\lambda$ and $A-YV-\lambda$
 turn out to be bijective in the end.

 Conversely, if the operators $A+V-\lambda$ and $A-YV-\lambda$ are in advance known to be bijective, then the graph subspace
 $\cG(\cH_0,X)$ is reducing for $A+V$ if and only if $Y$ is a strong solution to the operator Riccati equation.
\end{remark}

In Theorem \ref{thm:main}, certain regularity conditions on the operators $A+V$ and $A-YV$ are imposed. Since $A+V$ does not depend
on the choice of $Y$, the conditions on the operator $A+V$ are of an \emph{a priori} character, whereas the corresponding
conditions on $A-YV$ are of an \emph{a posteriori} character. For applications, the latter are rather inconvenient, so that one
naturally looks out for (stronger) a priori assumptions under which these a posteriori conditions are automatically satisfied. This
is discussed in the next section.

We close this section with the following remarks picking up the discussion at the end of Section \ref{sec:blockDiag}.

The first one addresses the issue of the first block diagonalization \eqref{eq:blockDiag1stDiss} in Theorem \ref{thm:blockDiag1st}.

\begin{remark}
 In accordance with Theorem \ref{thm:blockDiag1st}, there is a similar result to Theorem \ref{thm:main} addressing the first block
 diagonalization. In this case, one imposes regularity conditions on the operators $A+V$ and $A+VY$, and the property of $Y$ being
 a strong solution to the Riccati equation has to be replaced by the weaker notion from (ii) in Theorem \ref{thm:blockDiag1st}.
 However, if the perturbation $V$ is unbounded, there is no obvious additional assumption on $V$ to ensure the closedness of the
 operator $A+VY$. Hence, Remark \ref{rem:regCondClosed} can not be used to check the corresponding regularity conditions, which
 makes the result perhaps less suitable for applications. On the other hand, for the operator $A-YV$ closedness can be guaranteed
 in terms of relative bounds with respect to $A$, see Lemma \ref{lem:relBounds} and Theorem \ref{thm:mainRel} below.
\end{remark}

The second remark addresses the generalizations mentioned in Remarks \ref{rem:selfAdjoint} and \ref{rem:general}.

\begin{remark}
 The hypothesis that $\Dom(A)\subset\Dom(V)$ can be dropped in Theorem \ref{thm:main}, cf.\ also Remark \ref{rem:general}. However,
 in this case the technique described in Section \ref{sec:relBddPert} can not be applied directly to check the corresponding
 regularity conditions. The resulting statement is hence less suitable for our purposes.

 On the other hand, Theorem \ref{thm:main} admits a direct generalization to the case where the operator $A$ is not necessarily
 self-adjoint. Inspite of the reservations mentioned in Remark \ref{rem:selfAdjoint}, the corresponding result may still be of
 interest. We return to this more general case in Remark \ref{rem:moreGeneral}.
\end{remark}


\section{Relatively bounded perturbations}\label{sec:relBddPert}

As mentioned in the previous section, the hypotheses of Theorem \ref{thm:main} impose certain regularity conditions on the
operators $A+V$ and $A-YV$. Since in applications additional information on the operator $Y$ may not be available in advance, the
conditions on $A-YV$ tend to be hard to verify. In this section, we discuss stronger a priori assumptions on the perturbation $V$
that guarantee that these a posteriori type conditions on $A-YV$ are satisfied.

As a first step, we revisit the situation in \cite[Section 5]{AMM03} where the off-diagonal perturbation $V$ is assumed to be a
bounded operator.

\begin{theorem}{\cite[Theorem 5.3 and Theorem 5.5]{AMM03}}.\label{thm:AMM}
 Assume Hypotheses \ref{hyp1} and \ref{hyp2}. Furthermore, suppose that $V$ is bounded. Then, the graph subspace $\cG(\cH_0,X)$ is
 reducing for the operator $A+V$ if and only if $Y$ is a strong solution to the operator Riccati equation
 \[
  AY-YA-YVY+V=0\,.
 \]
 In this case, $A+V$  admits the block diagonalization
 \begin{equation}\label{eq:blockDiag1stBdd}
  T^{-1}(A+V)T = A+VY= \begin{pmatrix} A_0 + WX & 0\\ 0 & A_1-W^*X^*\end{pmatrix}\,.
 \end{equation}
 
 \begin{proof}
  Since $V$ and $Y$ are bounded, and $A$ is self-adjoint, the intersection of the resolvent sets of $A+V$ and $A-YV$ is not empty
  (see, e.g., \cite[Theorem V.4.3 and Problem V.4.8]{Kato66}). Hence, there exists some constant $\lambda$ such that $A+V-\lambda$
  and $A-YV-\lambda$ are bijective operators from $\Dom(A)$ to $\cH$. The stated equivalence now follows from Theorem
  \ref{thm:main}.

  Finally, the block diagonalization \eqref{eq:blockDiag1stBdd} is a direct consequence of equation \eqref{eq:blockDiag1st} in
  Theorem \ref{thm:blockDiag}.
 \end{proof}%
\end{theorem}

\begin{remark}\label{rem:AMM}
 The presented proof of Theorem \ref{thm:AMM} is based on Theorem \ref{thm:main}. A closer look at the proof of the latter shows
 that the property of the operator $T$ (resp.\ $T^*$) to map $\Dom(A)$ onto itself is essential for the equivalence in Theorem
 \ref{thm:AMM} to hold. In \cite[Lemma 5.3]{AMM03} this property has implicitly been used without justification.

 Moreover, the current proof is based on Lemma \ref{lem:invariance}, which is extracted from \cite[Theorem 5.5]{AMM03}. In other
 words, our proof of \cite[Lemma 5.3]{AMM03} requires some elements of \cite[Theorem 5.5]{AMM03}. In this sense, the two results
 \cite[Lemma 5.3]{AMM03} and \cite[Theorem 5.5]{AMM03} should not be considered as separate statements.

 In addition, the a posteriori character of the problem, which comes with the conditions on the operator $A-YV$ in Theorem
 \ref{thm:main}, is hidden in the original considerations in \cite[Section 5]{AMM03}. For bounded perturbations $V$, these
 conditions are automatically satisfied, which is why this a posteriori character is not visible in this case.
\end{remark}

With the preceding considerations in Sections \ref{sec:blockDiag} and \ref{sec:mainResult}, we now extend two aspects of Theorem
\ref{thm:AMM}. Not only do we present a second  block diagonalization for the operator $A+V$ besides \eqref{eq:blockDiag1stBdd}, we
also relax the boundedness requirement posed on $V$ (see Theorem \ref{thm:mainRel} below). Namely, taking into account Remark
\ref{rem:regCondClosed}, classic perturbation results like \cite[Theorem IV.3.17]{Kato66} allow to consider certain relatively
bounded perturbations $V$ with respect to $A$.

Recall that a linear operator $H\colon\cH\supset\Dom(H)\to\cH$ with $\Dom(H)\supset\Dom(A)$ is said to be \emph{$A$-bounded}
(or \emph{relatively bounded with respect to $A$}) if there exist constants $a,b\ge0$ such that
\begin{equation}\label{eq:defRelBound}
 \norm{Hx}\le a\norm{x} + b\norm{Ax} \quad\text{ for }\quad x\in\Dom(A)\,.
\end{equation}
If $H$ is $A$-bounded, the \emph{$A$-bound} of $H$ (or \emph{relative bound of $H$ with respect to $A$}) is defined as the infimum
of all possible choices for $b$ in \eqref{eq:defRelBound}, see \cite[Section IV.1.1]{Kato66}.

Note that in the situation of Hypotheses \ref{hyp1} and \ref{hyp2} the operator $V$ automatically is $A$-bounded, see
\cite[Remark IV.1.5]{Kato66} and \cite[Section V.3.3]{Kato66}.

The following result is based on Theorems IV.3.17, V.3.16, and V.4.3 from \cite{Kato66}.

\begin{lemma}\label{lem:relBounds}
 Assume Hypotheses \ref{hyp1} and \ref{hyp2}. If the operators $V$ and $YV$ are $A$-bounded with $A$-bound smaller than $1$, then
 $A+V$ is self-adjoint, $A-YV$ is closed, and the resolvent sets of $A+V$ and $A-YV$ are not disjoint, that is,
 $\rho(A+V)\cap\rho(A-YV)\neq\emptyset$.

 \begin{proof}
  Since $V$ is symmetric and $A$ is self-adjoint, it follows from \cite[Theorem V.4.3]{Kato66} that the operator $A+V$ also is
  self-adjoint if $V$ has $A$-bound smaller than $1$.

  Hence, it suffices to show the following statement: For every $A$-bounded linear operator $H\colon\cH\supset\Dom(H)\to\cH$ with
  $A$-bound smaller than $1$, the perturbed operator $A+H$ is closed, and there is a constant $k\ge0$ such that
  \[
   \ii\lambda \in \rho(A+H) \quad\text{ for }\quad \lambda\in\R\,,\ \abs{\lambda}>k\,.
  \]

  In order to prove this statement, choose $a\ge0$ and $0\le b < 1$ such that
  \[
   \norm{Hx} \le a\norm{x} + b\norm{Ax} \quad\text{ for }\quad x\in\Dom(A)\,.
  \]
  For $\lambda\neq 0$, it follows from \cite[Theorem V.3.16]{Kato66} that
  \begin{equation}\label{eq:resolventEstimates}
   \norm{(A-\ii\lambda)^{-1}} \le \frac{1}{\abs{\lambda}} \quad\text{ and }\quad \norm{A(A-\ii\lambda)^{-1}} \le 1\,.
  \end{equation}
  Define
  \[
   k:=\frac{a}{1-b}\ge0\,,
  \]
  and let $\lambda\in\R$ with $\abs{\lambda}>k$. Then, the estimates in \eqref{eq:resolventEstimates} imply that
  \begin{equation}\label{eq:relativeBoundEstimate}
   a\norm{(A-\ii\lambda)^{-1}} + b\norm{A(A-\ii\lambda)^{-1}} \le \frac{a}{\abs{\lambda}} + b < 1\,.
  \end{equation}
  Hence, by \cite[Theorem IV.3.17]{Kato66}, the operator $A+H$ is closed, and $\ii\lambda$ belongs to the resolvent set of $A+H$.
 \end{proof}%
\end{lemma}

\begin{remark}
 In the situation of Lemma \ref{lem:relBounds} observe that
 \[
  A-YV=A+V-TV\,.
 \]
 Hence, if $V$ has $A$-bound smaller than $1$, and therefore $A+V$ is self-adjoint, the condition that $YV$ has $A$-bound smaller
 than $1$ can be replaced by the condition that $TV$ has $(A+V)$-bound smaller than $1$.
\end{remark}

We are now ready to present the main result of this section. It is a direct extension of Theorem \ref{thm:AMM} from bounded to
relatively bounded perturbations $V$ with $A$-bound $0$.

\begin{theorem}\label{thm:mainRel}
 Assume Hypotheses \ref{hyp1} and \ref{hyp2}. Furthermore, suppose that $V$ is $A$-bounded with $A$-bound $0$. Then, the graph
 subspace $\cG(\cH_0,X)$ is reducing for the operator $A+V$ if and only if $Y$ is a strong solution to the operator Riccati
 equation
 \[
  AY-YA-YVY+V=0\,.
 \]
 In this case, $A+V$ is a self-adjoint operator that admits the two block diagonalizations
 \[
  T^{-1}(A+V)T = A+VY = \begin{pmatrix} A_0 + WX & 0\\ 0 & A_1-W^*X^*\end{pmatrix}
 \]
 and
 \[
  T^*(A+V)(T^*)^{-1} = A-YV = \begin{pmatrix} A_0 + X^*W^* & 0\\ 0 & A_1-XW\end{pmatrix}\,.
 \]
 In particular, the operators $A+VY$ and $A-YV$ are both closed and similar to the self-adjoint operator $A+V$. Moreover, the
 unitary block diagonalization \eqref{eq:unitaryBlockDiag} holds.

 \begin{proof}
  Since $V$ has $A$-bound $0$ and $Y$ is bounded, the operator $YV$ also has $A$-bound $0$. By Lemma \ref{lem:relBounds}, one has
  $\rho(A+V)\cap\rho(A-YV)\neq\emptyset$. Taking into account Remark \ref{rem:regCondClosed}, the statement then follows from
  Theorems \ref{thm:main}, \ref{thm:blockDiag}, and Remark \ref{rem:unitary}.
 \end{proof}%
\end{theorem}

\begin{remark}\label{rem:mainRel}
 Following Lemma \ref{lem:relBounds}, in Theorem \ref{thm:mainRel} it suffices to assume that both operators $V$ and $YV$ have
 $A$-bound smaller than $1$. Hence, if more information on the operator $Y$ is available, the hypothesis on the relative bound of
 $V$ can be weakened. If, for example, $Y$ is a contraction (see, e.g., \cite{ALT01} and \cite{KMM05} for a discussion on
 contractive solutions of the Riccati equation), the conclusion of Theorem \ref{thm:mainRel} holds if the relative bound of $V$ is
 smaller than $1$.
\end{remark}

\begin{example}\label{examp}
 Assume Hypothesis \ref{hyp1}. Suppose, in addition, that the spectra of $A_0$ and $A_1$ are subordinated, that is, 
 \[
  \sup \spec(A_0)<\inf \spec(A_1)\,.
 \]
 In particular, the operator $A_0$ is bounded from above and $A_1$ is bounded from below.

 If the operator $V$ is $A$-bounded with $A$-bound smaller than $1$, then the operator $A+V$ is self-adjoint. Thus, one can apply
 an extension of the Davis--Kahan $\tan2\Theta$ theorem (see \cite[Theorem 1]{MS06}) to conclude that 
 \[
  \Norm{\EE_A\bigl((-\infty, \sup \spec(A_0)]\bigr)-\EE_{A+V}\bigl((-\infty, \sup \spec(A_0)]\bigr)}\le\frac{\sqrt{2}}{2}\,,
 \]
 where $E_C(\sigma)$ stands for the spectral projection for a self-adjoint operator $C$ associated with a Borel set $\sigma$ on the
 real axis.
 
 In particular, the subspace $\Ran \EE_{A+V}\bigl((-\infty, \sup \spec(A_0)]\bigr)$ is a graph subspace with respect to
 $\cH_0:=\Ran \EE_A\bigl((-\infty, \sup \spec(A_0)]\bigr)$ (see, e.g., \cite[Corollary 3.4 (i)]{KMM03}), that is,
 \[
  \Ran \EE_{A+V}\bigl((-\infty, \sup \spec(A_0)]\bigr) = \cG(\cH_0, X)\,,
 \]
 where $X:\cH_0\to \cH_0^\perp$ is a linear contraction. Hence, the corresponding operator $Y$ defined as in Hypothesis \ref{hyp2}
 is a contraction as well. By Remark \ref{rem:mainRel}, the conclusion of Theorem \ref{thm:mainRel} holds providing a block
 diagonalization for the operator matrix $A+V$.
\end{example}

We close this work with the following generalizing observation.
\begin{remark}\label{rem:moreGeneral}
 In a formal sense, the self-adjointness of the operator $A$ is not required in Theorem \ref{thm:mainRel}, cf.\ also Remark
 \ref{rem:selfAdjoint}. Indeed, it suffices to assume that the operator $A$ is just closed and its resolvent set contains a half
 line $\mathfrak{h}$ such that the map $\mathfrak{h}\ni\mu\mapsto\abs{\mu}\norm{(A-\mu)^{-1}}$ is bounded. In this case,
 $\norm{(A-\mu)^{-1}}$ is arbitrarily small for $\mu\in\mathfrak{h}$ with large $\abs{\mu}$. Hence, with
 \[
  \norm{A(A-\mu)^{-1}}=\norm{(A-\mu)(A-\mu)^{-1}+\mu(A-\mu)^{-1}}\le 1 + \abs{\mu}\norm{(A-\mu)^{-1}}\,,
 \]
 the estimate \eqref{eq:relativeBoundEstimate} in the proof of Lemma \ref{lem:relBounds} can be replaced by
 \[
  a\norm{(A-\mu)^{-1}} + b\norm{A(A-\mu)^{-1}} \le a\norm{(A-\mu)^{-1}} + b\bigl(1 + \abs{\mu}\norm{(A-\mu)^{-1}}\bigr)<1
 \]
 for sufficiently small $b<1$ and $\mu\in\mathfrak{h}$ with sufficiently large absolute value.

 In this case, the block diagonalizations for $A+V$ remain valid, but the operator $A+V$ will not be self-adjoint in general.
\end{remark}


\end{document}